\documentstyle[12pt]{article}

\begin{document}
\author{S.V. L\"{u}dkovsky, J.C. Ferrando}
\title{Hahn-Banach theorems for $\kappa $-normed spaces}
\date{27 February 2002}
\maketitle

\begin{abstract}
For a new class of topological vector spaces, namely $\kappa $-normed
spaces, an associated quasisemilinear topological preordered space is
defined and investigated. This structure arise naturally from the
consideration of a $\kappa $-norm, that is a distance function between a
point and a $G_{\delta }$-subset. For it, analogs of the Hahn-Banach theorem
are proved.
\end{abstract}

\section{Introduction.}

Normable topological vector spaces play important role in functional
analysis, but their class does not encompass all locally convex spaces. A
new class of topological spaces, namely the class of $\kappa $-metric
topological spaces, was introduced earlier in works \cite{shep1,shep2},
where a regular $\kappa $-metric is defined as a non-negative function $\rho
:X\times 2_{o}^{X}\rightarrow \mathbf{R}$ satisfying axioms $(N1-4,5(b))$
below, $2_{o}^{X}$ being the family of all canonical closed subsets of $X$.
It is worth to mention that $\kappa $-metric spaces may be non-metrizable.
One of the most important examples of $\kappa $-metric spaces are locally
compact groups and generalized loop groups \cite{luumn48}. Different
distance functions for subsets of linear normed spaces were studied in \cite%
{att}. Topological vector spaces with $\kappa $-metrics satisfying
additional conditions related with the linearity of these spaces were
defined and studied in \cite{lusmz,ludkns}. For such $\kappa $-normed spaces
(see Definition 1.1 below) analogs of theorems about fixed point, closed
graph and open mapping have been proved. Free $\kappa $-normed spaces
generated by $\kappa $-metric uniform spaces with uniformly continuous $%
\kappa $-metrics have also been studied, as well as categorial properties of 
$\kappa $-normed spaces relative to products, projective and inductive
limits. Duality of $\kappa $-normed spaces and their applications has also
been investigated.

On the other hand, Hahn-Banach theorems play very important role for locally
convex spaces \cite{nari,schaef,topsoe}. In particular these theorems are
used for investigations concerning barrelledness properties of locally
convex spaces \cite{ferr}. In this work a structure of quasisemilinear
preordered topological space on $X\times S$, where $S$ is a family of
subsets of a topological vector space $X$, is naturally defined and
investigated. This structure in general does not reduce to a linear or even
a semilinear one. For it, analogs of the Hahn-Banach theorems are proved.
These theorems can be used for further investigations of analogs of
barrelledness for such quasisemiliniear spaces and their duality. Analogs of
separation theorems for quasisemilinear preordered spaces can be deduced
using the Hahn-Banach theorems. In turn, this can serve for a continuation
of investigations of applications of $\kappa $-normed spaces to differential
equations and economic problems (see also \cite{ludkns}).

\textbf{1.1. Definition.} A topological vector space $X$ over the filed $%
\mathbf{K}=\mathbf{R}$ or $\mathbf{C}$ or non-Archimedean supplied with the
family $S_{X}:=2_{o}^{X}$ of all canonical closed subsets or the family $%
S_{X}:=2_{\delta }^{X}$ of all closed $G_{\delta }$-subsets is called $%
\kappa $-normed if there exists on $X\times S_{X}$ a non-negative function $%
\rho (x,C)$, called a $\kappa $-norm, satisfying the following conditions:

$(N1)$ $\rho (x,C)=0$ if and only if $x\in C$;

$(N2)$ if $C\subset C^{\prime}$, then $\rho (x,C)\ge \rho (x,C^{\prime})$;

$(N3)$ the map $x\mapsto \rho (x,C)$ is uniformly continuous for each fixed $%
C\in S_{X}$;

$(N4)$ for each increasing transfinite sequence $\{C_{\alpha }\}$ with $%
C:=cl(\bigcup_{\alpha }C_{\alpha })\in S_{X}$ then $\rho (x,C)=\inf_{\alpha
}\rho (x,C_{\alpha })$, where $cl_{X}(A)=cl(A)$ denotes the closure of a
subset $A$ in $X$;

$(N5)(a)$ $\rho (x+y,cl(C_{1}+C_{2}))\leq \rho (x,C_{1})+\rho (y,C_{2})$ and

$(N5)(b)$ $\rho (x,C_{1})\leq \rho (x,C_{2})+{\bar{\rho}}(C_{2},C_{1})$
(with the maximum instead of the sum on the right sides of the inequalities
in the non-Archimedean case), where $\bar{\rho}(C_{2},C_{1}):=\sup_{x\in
C_{2}}\rho (x,C_{1})$;

$(N6)$ $\rho (\lambda x, \lambda C)= |\lambda |\rho (x,C)$ for each $\mathbf{%
K}\ni \lambda \ne 0$;

$(N7)$ $\rho (x+y,y+C)=\rho (x,C)$.

If we consider the empty set $\emptyset $ as an element of $S_{X}$, then

$(N8)$ $\rho (x,\emptyset )=\infty $ and $\rho (x,C)<\infty $ for each $x\in
X$ and $\emptyset \neq C\in S_{X}$.

The space $X\times S_{X}$ is called $\kappa $-normed if there is given a
fixed $\kappa $-norm $\rho $. We denote the $\kappa $-normed space by $%
(X,S_{X},\rho )$. For two closed subsets $A$ and $B$ in $X$ we denote by $A%
\hat{+}B$ the set $cl_{X}(A+B)$. If we put $k\circ A:=A\hat{+}...\hat{+}A$,
in general from $k\circ A=k\circ B$ does not follow $A=B$, and $kA$ may be
not equal to $k\circ A$ as may be easily seen. The main results of this
paper are Theorems 2.8, 2.10, 2.11 and 2.14.

\section{Hahn-Banach theorems for $\protect\kappa $-normed spaces.}

\textbf{2.1.1. Notes and Definitions.} Consider a family $S_{X}$ of subsets
of a topological vector space $X$ over the field $\mathbf{K}$, either $%
\mathbf{Q}$ or $\mathbf{R}$ or $\mathbf{C}$, such that

$(1)$ $\emptyset \in S_X$;

$(2)$ there is defined an addition `\thinspace $\hat{+}\,$' in $S_{X}$
making it a commutative semigroup with zero element, that is, $S_{X}^{2}\ni
(A,B)\mapsto A\hat{+}B\in S_{X}$ such that it is commutative and associative
with zero element $\emptyset $;

$(3)$ $S_{X}$ is preordered by the inclusion relation, that is, $A\leq B$ if
and only if $A\subset B$ or $A=B$, which obviously is $(3i)$ reflexive and $%
(3ii)$ transitive. Additionally we shall require that $(3iii)$ $A_{1}\leq
B_{1}$ and $A_{2}\leq B_{2}$ implies that $A_{1}\hat{+}B_{1}\leq A_{2}\hat{+}%
B_{2}$ and $(3iv)$ $A\hat{+}B\leq C\hat{+}B$ and $B\leq G$ implies that $A%
\hat{+}G\leq C\hat{+}G$;

$(4)$ it is invariant relative to the multiplicative group $\mathbf{K^{\ast }%
}:=\mathbf{K}\setminus \{0\}$, that is, $aA\in S_{X}$ for each $a\in \mathbf{%
K^{\ast }}$ and $A\in S_{X}$; moreover, it is semilinear relative to $%
\mathbf{K^{\ast }}$, that is, $a(A\hat{+}B)=aA\hat{+}aB$, $(ab)A=a(bA)$, $%
1A=A$ for each $A$ and $B$ in $S_{X}$ and each $a$ and $b$ in $\mathbf{%
K^{\ast }}$;

$(5)$ $S_X$ is invariant relative to shifts on all vectors from $X$, that
is, $v+A\in S_X$ for each $v\in X$ and $A\in S_X$.

Due to Condition $(4)$, $S_{X}$ is idempotent free, that is,

$(6)$ from $aA=aB$ for some $a\in \mathbf{K^*}$, it follows $A=B$.

Conditions $(3,4)$ imply that

$(7)$ $aA\leq aB$ for each $A\leq B$ in $S_{X}$ and $a\in \mathbf{K^{\ast }}$%
.

Conditions $(3,5)$ imply that

$(8)$ $A\le A\hat +B$ for each $A$ and $B$ in $S_X$ with $0\in B$.

If a family $S_{X}$ of subsets of $X$ satisfies Conditions $(1-5)$ we shall
call it a $\mathbf{K^{\ast }}$-quasisemilinear preordered space. Without $%
(3) $ we shall call it a $\mathbf{K^{\ast }}$-quasisemilinear space. If $A$%
\thinspace $\hat{+}B\leq A^{\prime }$\thinspace $\hat{+}B$ implies that $%
A\leq A^{\prime }$ we say that $S_{X}$ is monotonely cancellative.

If $S$ is an abstract family of sets of a topological space $X$ and there is
an abstract preorder `$\,\leq $\thinspace ' on $S$ such that the pair $(X,S)$%
, with the preorder $(x,A)\leq (y,B)$ if and only if $x=y$ and $A\leq B$,
may be equipped with an addition `\thinspace $+\,$' defined by $\left(
x,A\right) +\left( y,B\right) =\left( x+y,A\hat{+}B\right) $, where
`\thinspace $\hat{+}$\thinspace ' satisfies Conditions $(1-8)$, then $(X,S)$
will be called an (abstract) $\mathbf{K^{\ast }}$-quasisemilinear preordered
space. Providing $S$ with a topology such that the algebraic operations on $%
(X,S)$ are continuous; that is, relative to the topologies $\tau _{X}$ on $X$%
, $\tau _{S}$ on $S$ and $\tau _{\mathbf{K}}$ on $\mathbf{K}$, the maps $%
\hat{+}:(S,\tau _{S})^{2}\rightarrow (S,\tau _{S})$; $+:(X,\tau _{X})\times
(S,\tau _{S})\rightarrow (S,\tau _{S})$; $\mathbf{K^{\ast }}\times S\ni
(a,A)\mapsto aA\in S$ are continuous, then $(X,S)$ is called a topological $%
\mathbf{K^{\ast }}$-quasisemilinear space. If in $(1-8)$ $\mathbf{K^{\ast }}$
is substituted by $\mathbf{K_{+}}:=\mathbf{K}\cap (0,\infty )$, then it is
called a (topological) $\mathbf{K_{+}}$-quasisemilinear space.

Let $\mathbf{P}$ be either $\mathbf{K^{\ast }}$ or $\mathbf{K_{+}}$. For a
subset $J$ of a linear space $X$ over $\mathbf{K}$ we define, as usually, $%
sp_{\mathbf{P}}J:=\{x:x=\sum_{j=1}^{n}a_{j}p_{j};$ $a_{j}\in \mathbf{P},$ $%
p_{j}\in J,$ $n\in \mathbf{N}\}$. For a subset $J$ in $S_{X}$ we shall write 
$Sp_{\mathbf{P}}J:=\{A:A=a_{1}A_{1}\hat{+}a_{2}A_{2}\hat{+}...\hat{+}%
a_{n}A_{n};$ $a_{j}\in \mathbf{P},$ $A_{j}\in J,$ $n\in \mathbf{N}\}$.

A subset $F$ of $S_{X}$ is called hereditary if $A\in F$ whenever $A\leq B$
and $B\in F$. If $F\subset S_{X}$ and $S_{X}$ satisfies Conditions $(1-8)$,
then there exists a least $\mathbf{P}$-quasisemilinear preordered space
containing $F$, denoted by $\widehat{F}$, such that $\widehat{F}$ is the
intersection of all $\mathbf{P}$-quasisemilinear preordered subspaces of $%
S_{X}$ containing $F$. The set $\{A:A\in S_{X},$ $\mbox{there are}$ $B\in
S_{X}$ $\mbox{and}$ $C\in F$ $\mbox{such that}$ $A\hat{+}B\leq C\}$ is
called the hereditary face in $S_{X}$ generated by the subset $F$ of $S_{X}$.

For a linear space $X$ its family $\Omega _{X}$ of all singletons in $X$
satisfies Conditions $(1-8)$ with the ordinary sum and with the relation $%
x\leq y$ meaning $x=y$ in $X$. Hence, we get the notion of the hereditary
face of a subset $J$ in $X$. For $X\times S_{X}$ considered as $\Omega
_{X}\times S_{X}$ Conditions $(1-8)$ also are satisfied with the preorder $%
(x,A)\leq (y,B)$ if and only if $x=y$ and $A\leq B$, so we get the notion of
the hereditary face in $X\times S_{X}$ generated by a subset $P$ in $X\times
S_{X}$.

\textbf{2.1.2. Definitions and Conventions.} For a topological vector space $%
X$ and a family of subsets $S_{X}$ as in \S $2.1.1.$ a mapping $f:X\times
S_{X}\rightarrow \lbrack -\infty ,\infty ]$ is called sublinear if it
satisfies Conditions $(S1-S5)$:

$(S1)$ $f(0,\emptyset )=0$;

$(S2)$ $f(v+x,v+A)=f(x,A)$ for each $v$ and $x$ in $X$ and each $A\in S_X$;

$(S3)$ $f(x+y,\emptyset )\le f(x,\emptyset )+f(y,\emptyset )$ for each $x$
and $y$ in $X$;

$(S4)$ $f(0,A\hat +B)\le f(0,A)+f(0,B)$ for each $A$ and $B$ in $S_X$ with
the convention $\infty +(-\infty )=\infty $;

$(S5)$ $f(0,aA)=af(0,A)$ for each $a\in \mathbf{K_+}$.

It is called superlinear, if it satisfies Conditions $(S1,S2,S5)$ and

$(S3)^{\prime}$ $f(x+y,\emptyset )\ge f(x,\emptyset )+f(y,\emptyset )$ for
each $x$ and $y$ in $X$;

$(S4)^{\prime}$ $f(0,A\hat +B)\ge f(0,A)+f(0,B)$ for each $A$ and $B$ in $%
S_X $ with the convention $\infty +(-\infty )=-\infty $.

A mapping $f:X\times S_{X}\rightarrow \lbrack -\infty ,\infty ]$ is called
semilinear if it satisfies $(S1,S2,S5)$ and

$(S3)^{\prime \prime }$ $f(x,\emptyset )$ is linear by $x\in X$ on $X\times
\emptyset $;

$(S4)^{\prime \prime }$ $f(0,A\hat{+}B)=f(0,A)+f(0,B)$ for each $A$ and $B$
in $S_{X}$ and takes no more than one of the values $-\infty $, $\infty $.

Conditions $(S2,S5)$ imply that

$(S6)$ $f(ax,aA)=af(x,A)$ for each $a\in \mathbf{K_+}$, since $%
f(ax,aA)=f(0,a(A-x))=af(0,A-x)=af(x,A)$. From $(S2,S3)$ it follows that

$(S7)$ $f(x+y,A\hat +B)\le f(x,A)+f(y,B)$, since $f(x+y,A\hat
+B)=f(0,A-x\hat +B-y)\le f(0,A-x)+f(0,B-y)=f(x,A)+f(y,B)$.

In view of $2.1.1.(8)$ and $2.1.2.(S1-S4)$ or $2.1.2.(S1,S2,S3^{\prime
},S4^{\prime })$ for each nonvoid $A$ in $S$ there exists $x\in A$ such that 
$f(x,A)\geq 0$ for $f$ sublinear or semilinear, or $f(x,A)\leq 0$ for $f$
superlinear.

For an inequality $\sum a_i\le \sum b_j$ if the expression $\infty -\infty $
occurs on the left side we use the convention $\infty -\infty =-\infty $, if 
$\infty -\infty $ occurs on the right side we use the convention $\infty
-\infty =\infty $.

Let $F$ be a $\mathbf{P}$-quasisemilinear preordered subspace of $S_{X}$. A
mapping $f:F\rightarrow \lbrack -\infty ,\infty ]$ is called monotonely
cancellative if $f(A)\leq f(B)$ whenever $A\hat{+}C\leq B\hat{+}C$, with $A$
and $B$ in $F$, for some $C\in S_{X}$. A mapping $f:F\rightarrow \lbrack
-\infty ,\infty ]$ is called monotone if $f(A)\leq f(B)$ for each $A\leq B$
in $F$. These definitions extend trivially to mappings $f:H\times
F\rightarrow \left[ -\infty ,\infty \right] $ where $H$ is a linear subspace
of $X$ and $F$ is a set in $S_{X}$ satisfying Conditions \S $2.1.1.$ $(1-8)$
with $H$ instead of $X$.

\textbf{2.1.3. Definitions.} Let $p:X\times S_{X}\rightarrow (-\infty
,\infty ]$ be sublinear, $H$ be a linear subspace in $X$ and $F$ be a
subfamily in $S_{X}$ satisfying Conditions $(1-8)$ with $H$ instead of $X$
and let $g:H\times F\rightarrow (-\infty ,\infty ]$ be semilinear, monotone
and dominated by $p$. An extension $h$ of $g$ on $X\times S_{X}$, $h:X\times
S_{X}\rightarrow (-\infty ,\infty ]$, is called a monotone Hahn-Banach
extension if $h$ is semilinear, monotone and dominated by $p$.

If $p$ is sublinear, then from $p(x,A)<\infty $ and $p(y,B)<\infty $ it
follows $p(x+y,A\hat{+}B)\leq p(x,A)+p(y,B)<\infty $ and $%
p(ax,aA)=ap(x,A)<\infty $ for each $a\in \mathbf{K_{+}}$, hence the set $%
\{(x,A)\in X\times S_{X}:p(x,A)<\infty \}$ is a $\mathbf{P}$-quasisemilinear
preordered subspace of $X\times S_{X}$. Analogously, if $g$ is semilinear,
then the set $\{(x,A)\in X\times S_{X}:g(x,A)<\infty \}$ is a $\mathbf{P}$%
-quasisemilinear preordered subspace of $X\times S_{X}$.

\textbf{2.1.4. Remarks and Notations.} In \S $2.1.2.$ $(S5)$ the restriction 
$a>0$ is necessary since $f(0,aA)=af(0,A)$ for each $a\in \mathbf{K^{\ast }}$
would imply $f(0,A)=0$ for each $A$ in $S_{X}$ such that $A=-A$ (for
instance, if $A$ is a balanced subset $A$ in $X$).

If $\mathbf{P}$ is either $\mathbf{K^{\ast }}$ or $\mathbf{K_{+}}$, consider
the test relation for $x\in X$ and $A\in S_{X}$ consisting of the equality $%
y+nx=l+m$ and the inequality $B\hat{+}n_{1}\circ (a_{1}A)\hat{+}...\hat{+}%
n_{k}\circ (a_{k}A)\hat{+}C\leq J\hat{+}G\hat{+}C$ with $m$ in $X$, $y$ and $%
l$ in $H$, $B$ and $J$ in $F$, $\mathbf{Z}\ni n_{i}>0$ and $a_{i}\in \mathbf{%
K_{+}}$ with $i=1,...,k$, $k\in \mathbf{N}$, $n:=n_{1}a_{1}+...+n_{k}a_{k}$,
and $G$ and $C$ in $S_{X}$. A special test relation is characterized by the
condition $C=\emptyset $. Given as above a sublinear mapping $p:X\times
S_{X}\rightarrow (-\infty ,\infty ]$ and a semilinear and monotone mapping $%
g $ on $H\times F$ dominated by $p$, let us define the numbers $\xi (x,A)$
and $\xi _{0}(x,A)$ belonging to $[-\infty ,\infty ]$ by the following
formulas:

$\xi (x,A):=\inf \{ n^{-1}(p(m,G)+g(l,J)-g(y,B)):$ $y+nx=l+m$ $\mbox{and}$ $%
B\hat +n_1\circ (a_1A)\hat +...\hat + n_k\circ (a_kA)\hat +C\le J\hat +G\hat
+C \} $ and

$\xi _0(x,A):=\inf \{ n^{-1}(p(m,G)+g(l,J)-g(y,B)):$ $y+nx=l+m$ $\mbox{and}$ 
$B\hat + n_1\circ (a_1A)\hat +...\hat +n_k\circ (a_kA) \le J\hat +G \} $,

where the infimum is taken by all test relations and special test relations,
respectively. These functions are also denoted by $\xi _{p,g}(x,A)$ and $\xi
_{0,p,g}(x,A)$, respectively.

Let $E$ be $\mathbf{P}$-quasisemilinear subspace of $X\times S_X$. Denote by 
$\bigtriangleup E$ the set of all $x\in X$ and $A\in S_X$ for which there
exist $n_i\in \mathbf{N}$, $a_i\in \mathbf{K_+}$, $k\in \mathbf{N}$, $l$ and 
$y$ in $H:=\pi _X(E)$ with

$(i1)$ $y+nx=l$, \newline
and there exist $B$ and $C$ in $F:=\pi _{S}(E)$ and $Q\in S_{X}$ such that

$(i2)$ $B\hat +n_1\circ (a_1A)\hat +...\hat +n_k\circ (a_kA) \hat +Q=C\hat
+Q $, where variables $n$, $n_i$, $a_i$, $k$ are the same as above, $\pi _X:
X\times S_X\to X$ and $\pi _S: X\times S_X\to S_X$ are the natural
projections on $X$ and $S_X$ respectively. \newline
Also $\bigtriangleup _0E$ denotes the set of all $x\in X$ and all $A\in S_X$
for which there exist $n_i$, $a_i$, $k$, $m$ in $X$, $l$ and $y$ in $H$ with

$(ii1)$ $y+nx=l,$ \newline
and there exist $B$ and $C$ in $E$ such that

$(ii2)$ $B\hat +n_1\circ (a_1A)\hat +...\hat +n_k\circ (a_kA) =C$.

\textbf{2.1.5. Notes and Definitions.} Consider a $\kappa $-normed space $%
(X,S,\rho )$, where $X$ is a topological vector space with a topology $\tau
_{X}$, $S$ is a family of subsets such that each canonical closed subset
belongs to $S$ and $\rho :X\times S\rightarrow \lbrack 0,\infty )$ is a $%
\kappa $-norm. Let us define on $S$ several topologies.

Let $\tau _{1}$ denote the topology generated by the base $W(A,V):=\{B\in
S:A\subset B\subset A\hat{+}V\}$, where $A\in S$, $V\in \tau _{X}$, $0\in V$%
, $S=2_{o}^{X}$.

Let $\tau _{2}$ be the topology generated by the base $W(A,b):=\{B\in S:\bar{%
\rho}(A,B)<b\}$, where $\bar{\rho}(A,B):=\sup_{a\in A}\rho (a,B)$, $%
2_{o}^{X}\subset S\subset 2_{\delta }^{X}$.

The topology $\tau _{3}$ induced by the metric $D$ defined by the equation $%
D(A,B):=\bar{\rho}(A,B)+\bar{\rho}(B,A)$, $2_{o}^{X}\subset S\subset
2_{\delta }^{X}$.

For each $b>0$ and each $A\in S$ there exists a neighborhood $V$ of zero in $%
X$ such that $|\rho (x,A\hat{+}V)-\rho (x,A)|<b$ for each $x\in X$ (see
Corollary 5 \cite{lusmz}). If $A\subset B$ in $S$, then $\bar{\rho}(A,B)=0$,
hence $D(A,A\hat{+}V)<b$. Consequently, the topology $\tau _{1}$ on $S$ is
not weaker than $\tau _{3}$. Then consider on $X\times S$ the topologies $%
\zeta _{j}:=\tau _{X}\times \tau _{j}$, where $j=1,2,3.$

\textbf{2.1.6. Proposition.} \textit{Relative to each of the topologies }$%
\zeta _{j}$\textit{, where }$j=1,2,3,4$,\textit{\ }$X\times S$\textit{\ is a
topological }$\mathbf{P}$\textit{-quasisemilinear preordered space.}

\textbf{Proof.} Since the space $X\times S$ is a $\mathbf{P}$%
-quasisemilinear preordered space and $X$ is a topological vector space, it
remains to verify the continuity of algebraic operations in $X\times S$.

$\mathbf{\tau }$\textbf{$_{1}$.} For each $V\in \tau _{X}$ with $0\in V$
there are $V_{j}\in \tau _{X}$, $0\in V_{j}$, with $j=1$ and $j=2$ such that 
$V_{1}\hat{+}V_{2}\subset V$. If $A_{j}\subset B_{j}\subset A_{j}\hat{+}%
V_{j} $ for $j=1$ and $j=2$ with $A_{j}$ and $B_{j}$ in $S$, then $A_{1}\hat{%
+}A_{2}\subset B_{1}\hat{+}B_{2}\subset A_{1}\hat{+}A_{2}\hat{+}(V_{1}\hat{+}%
V_{2})$. Hence the addition $(A,B)\mapsto A\hat{+}B$ is continuous from $%
S^{2}$ to $S$. Besides, since $x+U\hat{+}W(A,V)\subset W(x+A,U+V)\subset
W(x+A,G)$, where $0\in U\cap V$, $U+V\subset G,$ $G\in \tau _{X}$, the
addition $(x,A)\mapsto x+A$ is continuous from $X\times S$ into $S$.

On the other hand, given $a\in \mathbf{K^{\ast }}$ and $A\in 2_{o}^{X}$, for
each $b>0$ and each neighbourhood $U$ of zero in $X$ there exists a
neighbourhood $V_{1}$ of zero in $X$ such that for each $c$ with $|a-c|<b$
then $cV_{1}\subset U$. If $A\subset B\subset A\hat{+}V_{1}$, then

$(i)$ $cA\subset cB\subset cA\hat{+}cV_{1}\subset cA\hat{+}U$. Since $X$ is
a topological vector space, then for each $0\in A\in 2_{o}^{X}$ and each $%
0\in V_{2}\in \tau _{X}$, $a\in \mathbf{K^{\ast }}$ and $b>0$ there exist $%
A_{1}\in 2_{o}^{X}$ and $0<\delta <b$ for which

$(ii)$ $aA\subset cA_{1}\subset aA\hat{+}V_{2}$ whenever $|a-c|<\delta $ and 
$c\neq 0$, moreover, $0\in A_{1}$. By the definition, if $B\in
cW(A_{1},V_{1})$ then there exists $B_{1}\in S$ such that $B=cB_{1}$ and $%
A_{1}\subset B_{1}\subset A_{1}\hat{+}V_{1}$, hence $cA_{1}\subset B\subset
cA_{1}\hat{+}V_{1}$. Choose $V_{1}$ and $V_{2}$ such that $V_{1}\hat{+}%
V_{2}\subset U$ and take $V=V_{1}\cap V_{2}$, then from Inclusions $(i,ii)$
it follows, that%
$$aA\subset cA_{1}\subset B\subset cA_{1}\hat{+}V_{1}\subset aA\hat{+}V_{1}%
\hat{+}V_{2}\subset aA\hat{+}U$$
for each $c$ such that $|a-c|<\delta $. Therefore, $\{c:|a-c|<\delta
\}W(A_{1},V)\subset W(aA,U)$ and hence the multiplication on scalars from $%
\mathbf{K^{\ast }}$ is continuous from $\mathbf{K^{\ast }}\times S$ into $S$.

$\mathbf{\tau }$\textbf{$_{2}$.} If $\bar{\rho}(A_{j},B_{j})<b/2$ for $A_{j}$
and $B_{j}$ in $S$, where $j=1,2$, then $\bar{\rho}(A_{1}\hat{+}A_{2},B_{1}%
\hat{+}B_{2})=\sup_{a_{1}\in A_{1},a_{2}\in A_{2}}\bar{\rho}%
(a_{1}+a_{2},B_{1}\hat{+}B_{2})\leq \sup_{a_{1}\in A_{1}}\bar{\rho}%
(a_{1},B_{1})+\sup_{a_{2}\in A_{2}}\bar{\rho}(a_{2},B_{2})=\bar{\rho}%
(A_{1},B_{1})+\bar{\rho}(A_{2},B_{2})<b$, hence $W(A_{1},b/2)\hat{+}%
W(A_{2},b/2)\subset W(A_{1}\hat{+}A_{2},b)$. From $\bar{\rho}(x+A,x+B)=\bar{%
\rho}(A,B)$ for each $x\in X$ and each $A$ and $B$ in $S$ we have $%
W(A,b)+x=W(x+A,b)$. In view of Lemma 2 \cite{lusmz} for each $A\in S$ there
exists a family $\{U_{j}:j\in \mathbf{N}\}$ such that $A=\bigcap_{j\in 
\mathbf{N}}U_{j}$, where $U_{j}\in \tau _{X}$ and $U_{j}\supset
cl_{X}(U_{j+1})$ for each $j$. Consider $A\in S$ such that $0\in A$. Then we
can take $0\in U_{j}$ for each $j$. Since $(X,\tau _{X})$ is the topological
vector space, denoting $A\hat{+}\left( -B\right) $ by $A\hat{-}B$ for each $j
$ there exists $V_{j}\in \tau _{X}$ such that $0\in V_{j}$, $V_{j}\hat{-}%
V_{j}\subset U_{j}$, $V_{j}\supset cl_{X}(V_{j+1})$ for each $j$. Therefore, 
$\bigcap_{j}V_{j}=:A_{1}\in S$ and $A_{1}\hat{-}A_{1}\subset A$. In view of 
\S $1.1.$ $(N3)$ for each such $A$ and each $b>0$ there exists $U$ such that 
$0\in U\in \tau _{X}$ and $|\rho (x+y,A_{1}\hat{-}A_{1})-\rho (x,A_{1}\hat{-}%
A_{1})|<b$ for each $y\in U$ and each $x\in X$. In view of Corollary 5 \cite%
{lusmz} this $U$ can be chosen such that $\bar{\rho}(A_{1}\hat{-}A_{1},B)<b$
for each $A_{1}\subset U$, $B\subset U$, $A_{1}\in S$, $0\in B\in \tau _{X}$%
. By \S $1.1.$ $(N5)$ we get the inequality%
$\bar{\rho}(y+A_{1},B)\leq \bar{\rho}(y+A_{1},A_{1}\hat{-}A_{1})+\bar{\rho}%
(A_{1}\hat{-}A_{1},B).$
Therefore, for each $0\in A\in S$ and each $b>0$ there exists $U\in \tau _{X}
$ such that $0\in U$ and there exists $0\in A_{1}\in S$ such that $%
x+y+W(A_{1},b)\subset x+W(A,2b)$ for each $x\in X$ and each $y\in U$. For
each $C\in S$ there exists $x\in C$ such that $C-x=:A\in S$ and $0\in S$.
Hence the addition $X\times S\ni (x,A)\mapsto x+A\in S$ is continuous. From $%
\bar{\rho}(ax,aB)=|a|\bar{\rho}(x,B)$ for each $a\in \mathbf{K^{\ast }}$ and
each $x\in X$ and each $B\in S$ we have $\bar{\rho}(aA,aB)=|a|\bar{\rho}(A,B)
$, hence $aW(A,b)=W(aA,|a|b)$ for each $a\in \mathbf{K^{\ast }}$ each $A\in S
$ and each $b>0$. It may be lightly seen that the multiplication is
continuous from $\mathbf{K^{\ast }}\times S$ into $S$.

$\mathbf{\tau }$\textbf{$_{3}$.} For each $b>0$ and $A\in S$ there exists $%
U\in \tau _{X}$ such that $\bar{\rho}(A,A\hat{+}U)<b$ (see Corollary 5 in %
\cite{lusmz}). Setting $W\left( A,b\right) =\left\{ B\in S:D\left(
A,B\right) <b\right\} $, since $D(x+A,x+B)=D(A,B)$ for each $x\in X$ and $%
A,B\in S$, then $x+W(A,b)=W(x+A,b)$ for each $x\in X$, each $A\in S$ and
each $b>0$. If $B\in S$ is such that $D(A,B)<b$, then $D(A+x,B)\leq
D(A+x,A)+D(A,B)<2b$ for each $x\in U$, since $\bar{\rho}(A\hat{+}U,A)=0$.
Hence the addition $(x,A)\mapsto x+A$ is continuous on $X\times S$. Since $%
D(A_{1}\hat{+}A_{2},B_{1}\hat{+}B_{2})\leq D(A_{1},B_{1})+D(A_{2},B_{2})$,
then the addition $(A,B)\mapsto A\hat{+}B$ is continuous from $S^{2}$ into $%
S $. From $\rho (ax,aA)=|a|\rho (x,A)$ for each $a\in \mathbf{K^{\ast }}$
and $x\in X$ and $A\in S$, it follows that $D(aA,aB)=|a|D(A,B)$ for each $A$
and $B$ in $S$ and each $a\in \mathbf{K^{\ast }}$. Therefore, the
multiplication on scalars $(a,A)\mapsto aA$ is continuous from $\mathbf{%
K^{\ast }}\times S$ into $S$.

\textbf{2.2. Lemma.} $(I)$. $\bigtriangleup E$ \textit{and }$\bigtriangleup
_{0}E$\textit{\ are }$\mathbf{P}$\textit{-quasisemilinear subspaces such
that }$\bigtriangleup (\bigtriangleup E)=\bigtriangleup E$\textit{, }$%
\bigtriangleup _{0}(\bigtriangleup _{0}E)=\bigtriangleup _{0}E$\textit{. }

$(II)$.\textit{\ If} $g$ \textit{is a finite semilinear and monotonely
cancellative map on }$E$\textit{, then there exists a unique semilinear and
monotonely cancellative extension} $h$ \textit{on} $\bigtriangleup E$\textit{%
. }

$(III)$\textit{. If }$g$\textit{\ is a finite semilinear and monotone map
defined on }$E$\textit{, then }$g$ \textit{has a unique semilinear and
monotone extension to} $\bigtriangleup _{0}E$\textit{.}

\textbf{Proof.} $(I)$. Equalities $2.1.4.(i,ii)$ are invariant relative to
the multiplication on scalars from $\mathbf{P}$ and also relative to shifts
on vectors $x$ from $X$, hence $\bigtriangleup E$ and $\bigtriangleup _{0}E$
satisfy Conditions $2.1.1.(4,5)$. Addition of two Equalities of type $%
2.1.4.(i1,i2)$ for $y_{j},x_{j},l_{j}$ and also for $B_{j}$, $A_{j}$, $C_{j}$
and $Q_{j}$ with $j=1$ and $j=2$ shows that $\bigtriangleup E$ satisfies
Condition $2.1.1.(2)$; analogously for $\bigtriangleup _{0}E$ with
Equalities $2.1.4.(ii1,ii2)$. From this Statement $(I)$ follows.

$(II)$. Let $g$ be defined on $E$ and $(x,A)$ be satisfying Conditions $%
2.1.4.(i1,i2)$. Then $g(y,B)$ and $g(l,J)$ are in $\mathbf{R}$. Therefore,
defining $g(x,A)$ by $g(x,A)=n^{-1}(g(l,C)-g(y,B))$ and using semilinearity
of $g$ on $\mathbf{P}$-quasisemilinear $E$ we get $g(x,A)$ defined on $(sp_{%
\mathbf{P}}x)\times (Sp_{\mathbf{P}}A)$. A construction of $g$ on $%
\bigtriangleup E$ can be done with the help of Conditions $2.1.4.(i1,i2)$ by
considering all $(x_{k},A_{k})$ in $\bigtriangleup E\setminus E$ such that $%
g $ on $\hat{E}_{k}$ with $\hat{E}_{k}:=E\cup (x_{k},A_{k})$ is defined by
assigning $g_{k}=g(x_{k},A_{k})\in (-\infty ,\infty ]$ such that $g_{k}\leq
g_{l}$ for each $(x_{k},A_{k})\leq (x_{l},A_{l})$ in $\bigtriangleup E$.
Considering the family $\Upsilon :=\{(\Sigma ,h):h|_{E}=g\}$ of $\mathbf{P}$%
-quasisemilinear subspaces $\Sigma $ in $\bigtriangleup E$ and semilinear
and monotonely cancellative mappings $h$ on $\Sigma $ ordered by $(\Sigma
_{1},h_{1})\leq (\Sigma _{2},h_{2})$ if and only if $\Sigma _{1}\subset
\Sigma _{2}$ and $h_{2}|_{\Sigma _{1}}=h_{1}$, due to the Kuratowski-Zorn
lemma we get that there exists a maximal element $(\Sigma ,h)\in \Upsilon .$
If for this maximal element $(\Sigma ,h)$ would be $\Sigma \neq
\bigtriangleup E$, then the construction with $(x_{k},A_{k})\in
\bigtriangleup E\setminus \Sigma $ could be continued contradicting $(\Sigma
,h)$ maximality, hence there exists $(\Sigma ,h)\in \Upsilon $ such that $%
\Sigma =\Delta E$ and $h|_{E}=g$.

$(III)$. The proof of the last statement is analogous to that of $(II)$ with
the help of Condition $2.1.4.(ii)$ instead of $2.1.4.(i)$.

\textbf{2.3. Lemma.} \textit{Let }$M$\textit{\ be the hereditary face
generated by }$P$\textit{, where }$P:=\{(x,A)\in X\times S_{X}:g(x,A)<\infty
\}\cup \{(x,A)\in X\times S_{X}:p(x,A)<\infty \}$. \textit{Then there exists
a monotone Hahn-Banach extension of }$g$\textit{\ on }$X\times S_{X}$\textit{%
\ if and only if }$M\cap \{(x,A):g(x,A)=\infty \}=\emptyset ,$\textit{\ and
the restriction of }$g$\textit{\ to }$M\cap E$\textit{, where }$E=H\times F$%
\textit{\ (see \S 2.1.3) has a (necessary finite) semilinear and monotone
extension on }$M$\textit{\ dominated by }$p|_{M}$.

\textbf{Proof.} If $Q:=M\cap \{(x,A):g(x,A)=\infty \}\neq \emptyset $, then
evidently for the extension $g$ on $X\times S_{X}$ the condition $%
range(g)\subset (-\infty ,\infty )$ would not be satisfied. On the other
hand, if $Q=\emptyset $, such extension exists as in Lemma 2.2.

\textbf{2.4. Lemma.} \textit{If }$h$\textit{\ is a monotone Hahn-Banach
extension of }$g$\textit{, then }$h\leq \xi _{0}$\textit{. If additionally }$%
h$\textit{\ is finite or, more generally, if }$h$\textit{\ is monotonely
cancellative, then} $h\leq \xi $.

\textbf{Proof.} If $B\hat{+}n_{1}\circ (a_{1}A)\hat{+}...\hat{+}n_{k}\circ
(a_{k}A)\leq E\hat{+}G$ and $y+nx=l+m$ with $n=n_{1}a_{1}+\ldots +n_{k}a_{k}$%
, then $h(y+nx,B\hat{+}n_{1}\circ (a_{1}A)\hat{+}...\hat{+}n_{k}\circ
(a_{k}A))=h(y,B)+nh(x,A)$ $\leq h(l+m,E\hat{+}G)=h(l,E)+h(m,G)$. Therefore,

$(i)$ $h(x,A)\leq n^{-1}\{h(l,E)+h(m,G)-h(y,B)\}$ $=n^{-1}%
\{h(m,G)+g(l,E)-g(y,B)\}$ $\leq n^{-1}\{p(m,G)+g(l,E)-g(y,B)\}$ for each
special test relation, consequently, $h(x,A)\leq \xi _{0}(x,A)$ for each $%
(x,A)\in X\times S_{X}$. If in addition $h$ is monotonely cancellative and $B%
\hat{+}n_{1}\circ (a_{1}A)\hat{+}...\hat{+}n_{k}\circ (a_{k}A)\hat{+}C\leq E%
\hat{+}G\hat{+}C$ and $y+nx=l+m$, then $h(y+nx,B\hat{+}n_{1}\circ (a_{1}A)%
\hat{+}...\hat{+}n_{k}\circ (a_{k}A)\hat{+}C)\leq h(l+m,E\hat{+}G\hat{+}C)$
implies $h(y+nx,B\hat{+}n_{1}\circ (a_{1}A)\hat{+}...\hat{+}n_{k}\circ
(a_{k}A))\leq h(l+m,E\hat{+}G)$. So Inequalities $(i)$ for each test
relation give $h(x,A)\leq \xi (x,A)$ for each $(x,A)\in X\times S_{X}$.

If $h$ is finite, then $h\left( x,C\right) \in \mathbf{R}$ for each $C$,
consequently, if $A\hat{+}C\leq B\hat{+}C$, then $h(x,A)+h(x,C)\leq
h(x,B)+h(x,C)$ and $h(x,A)\leq h(x,B)$, hence $h$ is monotonely cancellative.

\textbf{2.5. Lemma.} \textit{If} $\xi >-\infty $ \textit{(see \S 2.1.4),
then }

$(i)$ $\xi $ \textit{is sublinear and monotonely cancellative; }

$(ii)$ $\xi \leq p$\textit{; }

$(iii)$ $g$ \textit{is monotonely cancellative; }

$(iv)$ $\xi _{p,h}=\xi _{p,g}$\textit{; }

$(v)$ $\xi _{\xi ,g}=\xi $.\textit{\ }

\textit{In the case of} $\xi _{0}$ \textit{and }$\bigtriangleup _{0}E$ 
\textit{instead of }$\xi $\textit{\ and }$\bigtriangleup E$\textit{\ the
same statements hold, but in }$(i,iii)$\textit{\ only monotonicity is
guaranteed.}

\textbf{Proof.} $(i)$. Since $g$ is semilinear and $p$ is sublinear, then $%
\xi (ax,aA)=a\xi (x,A)$ for each $a\in \mathbf{K_{+}}$ and $(x,A)\in X\times
S$. Evidently, $\xi $ also satisfies Conditions $(S1,S2)$. Since $p$
satisfies $(S3,S4)$ and $g$ satisfies $(S3,S4)^{\prime \prime }$, then $\xi $
satisfies $(S3,S4)$. To prove monotone cancellation property of $\xi $ for
each test relation $(y+nx,B\hat{+}n_{1}\circ (a_{1}A)\hat{+}...\hat{+}%
n_{k}\circ (a_{k}A)\hat{+}C)\leq (l+m,J\hat{+}G\hat{+}C)$ with $A\hat{+}%
D\leq A^{\prime }\hat{+}D$ for some $D\in S$ it is sufficient to find a test
relation $(y+nx,B\hat{+}n_{1}\circ (a_{1}A^{\prime })\hat{+}...\hat{+}%
n_{k}\circ (a_{k}A^{\prime })\hat{+}C^{\prime })\leq (l+m,J\hat{+}G\hat{+}%
C^{\prime })$. Let $C\in S$, put $C\hat{-}C:=C\hat{+}(-1\cdot C)$ as usual
and note that $-1\cdot \left( C\hat{-}C\right) =C\hat{-}C$ and $0\in C\hat{-}%
C$. Since $S$ is a conmutative semigroup, then $C=A\hat{+}D$ for $A$ and $D$
in $S$. So, due to Conditions \S $2.1.1.(2,8)$ we get $A^{\prime }\hat{+}%
D\leq (A^{\prime }\hat{+}D)\hat{+}[(A\hat{+}D)\hat{-}(A\hat{+}D)]=(A\hat{+}D)%
\hat{+}(A^{\prime }\hat{-}A)\hat{+}(D\hat{-}D)$, so it suffices to take $%
C^{\prime }=C\hat{+}n_{1}\circ (a_{1}Q)\hat{+}...\hat{+}n_{k}\circ (a_{k}Q)$
with $Q:=[(A^{\prime }\hat{-}A)\hat{+}(D\hat{-}D)].$

$(ii)$. Take a test relation such that $l=0$, $y=0$ and $m=nx$, $E=\emptyset 
$, $B=\emptyset $, $G\ge n_1\circ (a_1A)\hat +...\hat +n_k\circ (a_kA)$,
then $\xi (x,A)\le n^{-1}p(m,G) \le p(x,A)$, consequently, $\xi \le p$.

$(iii)$. By the conditions of \S $2.1.3.$ $g$ is monotone and $g\leq p$,
hence, by Lemma 2.4 $g$ is monotonely cancellative.

$(iv)$. $\xi _{p,g}$ is defined on $\bigtriangleup E$, since each test
relation is a particular case of conditions defining $\bigtriangleup E$.
Each test relation with $g$ and $E$ is certainly a test relation with $h$
and $\bigtriangleup E$, consequently, $\xi _{p,h}\le \xi _{p,g}$. Consider
Conditions $2.1.4.(i1,i2)$: $y^{\prime}+nx=l^{\prime}$ and $B^{\prime}\hat
+k\circ A\hat +Q=C^{\prime}\hat +Q$, then the test relation of \S 2.1.4
takes the form $ky+ky^{\prime}+knx=kl+kl^{\prime}+km$ and $k\circ B\hat
+k\circ B^{\prime}\hat +k\circ ( n_1\circ (a_1A)\hat +...\hat +n_k\circ
(a_kA))\hat +k\circ C \le k\circ J\hat +k\circ C^{\prime}\hat +k\circ G\hat
+k\circ C$. But this shows that each test relation with $h$ and $%
\bigtriangleup E$ is also a test relation with $g$ and $E$, since $g$ is
semilinear, hence $\xi _{p,h}=\xi _{p,g}$.

$(v)$. From $(ii,iv)$ we get $\xi _{\xi ,g}\le \xi _{p,g}$. On the other
hand, composition of two subsequent test relations is also a test relation
and from the definition of $\xi $ this statement follows.

For $\xi _0$ and $\bigtriangleup _0E$ instead of $\xi $ and $\bigtriangleup
E $ only monotonicity in $(i,iii)$ is guaranteed, since for $\xi _0$ may be
a term $p(m,G)=\infty $, but for $\xi $ due to variation of $C$ this term
can be chosen $p(m,G)<\infty $, otherwise $p=\infty $ and $\xi =\infty $ for
all arguments that also gives monotone cancellation property of $\xi $.

\textbf{2.6. Lemma.} \textit{Let }$\eta $\textit{\ be finite superlinear and
dominated by }$\xi $\textit{. For each }$v\in X$ \textit{and each }$D\in S$%
\textit{\ define }

$u_{v,D}:=\sup \{n^{-1}(\eta (x_{1},A_{1})-\xi
(m_{1},G_{1})-g(l_{1},J_{1})+g(y_{1},B_{1})):$ $y_{1}+x_{1}=l_{1}+m_{1}+nv,$ 
$B_{1}\hat{+}A_{1}\hat{+}C_{1}\leq J_{1}\hat{+}G_{1}\hat{+}n_{1}\circ
(a_{1}D)\hat{+}...\hat{+}n_{k}(a_{k}D)\hat{+}C_{1}\}$;

$U_{v,D}:=\inf \{n^{-1}(-\eta (x_{2},A_{2})+\xi
(m_{2},G_{2})-g(l_{2},J_{2})-g(y_{2},B_{2})):$ $y_{2}+x_{2}+nv=l_{2}+m_{2}$%
\textit{,} $B_{2}\hat{+}A_{2}\hat{+}n_{1}\circ (a_{1}D)\hat{+}...\hat{+}%
n_{k}\circ (a_{k}D)$ $\hat{+}C_{2}\leq J_{2}\hat{+}G_{2}\hat{+}C_{2}\}$%
\textit{; where }$(y_{j},B_{j})$\textit{\ and }$(l_{j},J_{j})$\textit{\ with}
$j=1$\textit{\ and }$j=2$ \textit{are in }$E$\textit{, }$(x_{j},A_{j})$%
\textit{, }$(m_{j},G_{j})$\textit{\ and }$(z_{j},C_{j})$\textit{\ with }$j=1$%
\textit{\ and }$j=2$\textit{\ are in }$X\times S$\textit{. Then }

$(i)$ $\eta (m+v,D\hat{+}G)-\xi (m,G)\leq u_{v,D}\leq U_{v,D}\leq \xi (x+v,A%
\hat{+}D)-\eta (x,A)$\textit{\ for each }$(x,A)$\textit{\ and }$(m,G)$%
\textit{\ in }$X\times S$\textit{; }

$(ii)$ $-\infty <u_{v,D}\leq U_{v,D}\leq \infty $ \textit{and }$%
U_{v,D}<\infty $\textit{\ if }$(v,D)$\textit{\ lies in the hereditary face
generated by }$\{(z,P):\xi (z,P)<\infty \}$\textit{; }

$(iii)$ $u_{V,D}$ \textit{and }$U_{v,D}$\textit{\ are the same for }$p$%
\textit{\ instead of }$\xi $\textit{.}

\textbf{Proof.} Taking $y=0$, $l=0$, $J=\emptyset $, $B=\emptyset $ and
omitting indices from the definition of $U_{v,D}$ we get $U_{v,D}\le
n^{-1}(-\eta (x,A)+\xi (m,G)))$. Then for $n=1$, $m=x+v$ and $A\hat +D\hat
+C=G\hat +C$ due to Lemma $2.5.(i)$ $\xi (m,G)=\xi (x+v,A\hat +D)$, hence $%
U_{v,D}\le \xi (x+v,A\hat +D)-\eta (x,A)$.

Choose now $l_1=0$, $y_1=0$, $n=1$, $B=\emptyset $, $J=\emptyset $, $x=m+v$, 
$A=G\hat +D$, then $\eta (x,A)=\eta (m+v,D\hat +G)$ and from the definition
of $u_{v,D}$ we get $\eta (m+v,D\hat +G)-\xi (m,G)\le u_{v,D}$.

From the equalities $y_{1}+x_{1}=l_{1}+m_{1}+nv$ and $%
y_{2}+x_{2}+nv=l_{2}+m_{2}$ it follows $%
y_{1}+x_{1}-l_{1}-m_{1}=l_{2}+m_{2}-y_{2}-x_{2}$. From the inequalities $%
B_{1}\hat{+}A_{1}\hat{+}C_{1}\leq J_{1}\hat{+}G_{1}\hat{+}n_{1}\circ (a_{1}D)%
\hat{+}...\hat{+}n_{k}\circ (a_{k}D)\hat{+}C_{1}$ and $B_{2}\hat{+}A_{2}\hat{%
+}n_{1}\circ (a_{1}D)\hat{+}...\hat{+}n_{k}\circ (a_{k}D)\hat{+}C_{2}\leq
J_{2}\hat{+}G_{2}\hat{+}C_{2}$ it follows $B_{1}\hat{+}A_{1}\hat{+}C_{1}\hat{%
+}B_{2}\hat{+}A_{2}\hat{+}C_{2}\leq $ $J_{1}\hat{+}G_{1}\hat{+}n_{1}\circ
(a_{1}D)\hat{+}...\hat{+}n_{k}\circ (a_{k}D)\hat{+}C_{1}\hat{+}B_{2}\hat{+}%
A_{2}\hat{+}C_{2}$ $\leq $ $J_{2}\hat{+}G_{2}\hat{+}C_{2}\hat{+}J_{1}\hat{+}%
G_{1}\hat{+}C_{1}$. Taking $C:=C_{1}\hat{+}B_{2}\hat{+}A_{2}\hat{+}C_{2}$
for the first inequality and $C:=C_{1}\hat{+}C_{2}\hat{+}J_{1}\hat{+}G_{1}$
for the second inequality we get $u_{v,D}\leq U_{v,D}$.

$(ii)$. The inequality $u_{v,D}>-\infty $ follows from finiteness of $\eta $
and $g$ and taking $m_1=0$, $G_1=\emptyset $ for which $\xi (0,\emptyset )=0$%
. Then $U_{v,D}<\infty $ if there exists $\xi (m,G)<\infty $, that is the
case in the hereditary face.

$(iii)$. Using test relations it is possible to take for each $b>0$ a test
relation such that $|\xi (m,G)-g(l,J)+g(y,B)-p(m,G)|<b$. Using semilinearity
of $g$ and substituting arguments $B^{\prime}=B\hat +B_j$, $J^{\prime}=J\hat
+J_j$, $l^{\prime}=l+l_j$, $y^{\prime}=y+y_j$ we get that $u_{v,D}$ and $%
U_{v,D}$ are the same for $p$ instead of $\xi $.

\textbf{2.7. Lemma.} \textit{Assume that }$\eta $\textit{\ is finite
superlinear mapping on }$X\times S$\textit{\ and such that }$\xi \geq g\geq
\eta $ \textit{on }$E$\textit{, where }$g$\textit{\ and} $\xi =\xi _{p,g}$%
\textit{\ are defined on }$E$\textit{. Let }$(v,D)\notin \bigtriangleup E$ 
\textit{and }$F:=\{(z,P):z=v_{1}+e,$ $P=D_{1}\hat{+}C$\textit{$,$ }$(e,C)\in
E,$ $e\in X,$ $C\in S,$ $v_{1}\in sp_{\mathbf{K}}v,$ $D_{1}\in Sp_{\mathbf{P}%
}D\}$\textit{. For each }$\gamma \in \mathbf{R}$\textit{\ denote by }$%
h_{\gamma }$\textit{\ the uniquely determined semilinear extension of }$g$%
\textit{\ from }$E$\textit{\ to }$F$\textit{\ such that }$h_{\gamma
}(v,D)=\gamma $\textit{. Then }

$(i)$ $\xi _{p,h_{\gamma }}\geq h_{\gamma }\geq \eta $\textit{\ on }$F$%
\textit{\ if and only if }

$(ii)$ $u_{v,D}\leq \gamma \leq U_{v,D}$\textit{.}

\textbf{Proof.} The existence of $h_{\gamma }$ follows from Lemma 2.2. By
Lemma 2.5 $\xi _{p,h_{\gamma }}\ge \eta $ on $E.$ Suppose $(i)$ is satisfied
on $F$ and $(v,D)\notin \bigtriangleup E$. In view of Lemma 2.5 $h_{\gamma }$
is monotonely cancellative. Then $h_{\gamma }(m_2,G_2)+h_{\gamma
}(l_2,J_2)-h_{\gamma } (y_2,B_2)$ $\ge h_{\gamma }(x_2+nv,A_2\hat + n_1\circ
(a_1D) \hat +...\hat +n_k\circ (a_kD))=h_{\gamma }(x_2,A_2)+n\gamma $, since 
$h$ is semilinear. On the other hand, $-\eta (x_2,A_2) \ge -\xi (x_2,A_2)$,
consequently, $U_{v,D}\ge \gamma $.

Using $\xi \ge \eta $, $h_{\gamma }\ge \eta $ and monotone cancellation of $%
h_{\gamma }=:h$ we get $\eta (x_1,A_1)-\xi _{p,h}(m_1,G_1)-$ $%
h(l_1,J_1)+h(y_1,B_1)\le $ $\eta (x_1,A_1)-h(m_1,G_1)-h(l_1,J_1)+ h(y_1,B_1)$
$\le h(x_1,A_1)+h(y_1,B_1)-h(m_1,G_1)-h(l_1,J_1)$ $=h(x_1+y_1,A_1\hat
+B_1)-h(m_1+l_1,G_1\hat +J_1)$ $\le nh(v,D)=n\gamma $, since $%
h(x_1,A_1)+h(y_1,B_1)=$ $h(x_1+y_1,A_1\hat +B_1)\le h(l_1+m_1+nv, J_1\hat
+G_1\hat +n_1(a_1D)\hat +...\hat +n_k(a_kD))=$ $%
h(l_1,J_1)+h(m_1,G_1)+nh(v,D) $, consequently $u_{v,D}\le \gamma $.

Suppose now that $(ii)$ is satisfied and $\xi \ge g\ge \eta $ on $E$. In
view of Lemma 2.5 we have $\xi _{p,h}\ge h$ on $\bigtriangleup E$. Applying
Lemma 2.4 to $-g$ and $-\eta $ we have $h\ge \eta $ on $E$. It remains to
consider $F\ne \bigtriangleup E$. In view of Lemma $2.6.(i)$ we have $\eta
(v,D)\le u_{v,D}$ (taking $m=0$ and $G=\emptyset $) also $\xi (v,D)\ge
U_{v,D}$ (taking $x=0$ and $A=\emptyset $). Combining this with $(ii)$ we
get $\xi (v,D)\ge \gamma \ge \eta (v,D)$.

\textbf{2.8. Theorem. } \textit{Let }$p:X\times S_{X}\rightarrow \lbrack
-\infty ,\infty ]$\textit{\ be sublinear and }$g:H\times F\rightarrow
\lbrack -\infty ,\infty ]$\textit{\ be a finite semilinear, monotone and
dominated by }$p$\textit{\ mapping as above. }

$(i)$\textit{. If $X\times S$ is the hereditary face generated by }$H\times F\cup
\{(x,C):p(x,C)<\infty \}$\textit{\ or more generally by }$\{(x,C):\xi
(x,C)<\infty \}$\textit{, then a necessary and sufficient condition that }$g$%
\textit{\ has a finite Hahn-Banach extension }$h$\textit{\ on }$X\times
S_{X} $\textit{\ is that there exists a finite superlinear map }$\eta $%
\textit{\ on }$X\times S_{X}$\textit{\ such that }$\eta \leq \xi $.\textit{\
If this condition is fulfilled, then }$h$\textit{\ can be chosen such that }$%
\eta \leq h\leq \xi $\textit{. }

$(ii)$\textit{. Then a necessary and sufficient condition that }$g$\textit{\
has a monotone Hahn-Banach extension $h$ on }$X\times S_{X}$\textit{\ is
that, there exists a superlinear map} $\eta :X\times S_{X}\rightarrow
(-\infty ,\infty ]$ \textit{such that }$j\leq \xi _{0}$\textit{. When this
condition is fulfilled, then }$h$\textit{\ can be chosen such that }$\eta
\leq h\leq \xi _{0}$\textit{\ on }$X\times S_{X}$\textit{. }

$(iii)$\textit{. If the preordering on }$X\times S$\textit{\ is equality,
then }$g$\textit{\ has a finite Hahn-Banach extension on $X\times S$ if and
only if there exists a superlinear map }$\eta :X\times S\rightarrow (-\infty
,\infty ]$\textit{\ such that }$\eta \leq \xi _{0}$\textit{\ and }$g$\textit{%
\ is cancellative.}

\textbf{Proof.} $(i)$. Necessity follows from Lemma 2.4. We prove
sufficiency. Consider a maximal pair $(F,h)$ with $F$ a $\mathbf{P}$%
-quasisemilinear subspace of $\mathbf{P}$-quasisemilinear preordered space $%
X\times S$ such that $F\supset E$ and $h$ is a finite semilinear extension
of $g$ with $\xi _{p,h}\geq \eta $. By Lemma $2.5.(iii)$ $h$ is monotone on $%
F$. By Lemma $2.5.(v)$ $\bigtriangleup F=F$. Due to Lemmas $2.6$ and $2.7$ $%
\bigtriangleup F=X\times S$. Here is employed the assumption that $X\times S$
is the hereditary face generated by $\{\xi <\infty \}$, which gives $%
[u_{v,D},U_{v,D}]\neq \emptyset .$ By Lemmas $2.4$ and $2.7$ $\eta \leq
h\leq \xi $ on $X\times S$.

$(ii)$. Necessity is obvious. We prove sufficiency. Assume $\xi _{0}\geq
\eta $ with a superlinear map $\eta :X\times S_{X}\rightarrow (-\infty
,\infty ]$. Consider the face $Z=Z_{1}\times Z_{2}$ in $X\times S_{X}$
generated by $\{(x,A):\xi _{0}(x,A)<\infty \}$, where $Z_{1}\subset X$ and $%
Z_{2}\subset S_{X}$. For this $Z$ put $p^{\prime }:=p|_{Z}$. Let $\xi
^{\prime }$ be defined for such $p^{\prime }$ and $g$ as above. We will show
that $\xi ^{\prime }$ dominates the restriction $\eta |_{Z}$. For this it is
necessary to prove that when $B\hat{+}n_{1}\circ (a_{1}A)\hat{+}...\hat{+}%
n_{k}\circ (a_{k}A)\hat{+}C\leq E\hat{+}G\hat{+}C$ and $y+nx=l+m$ is
satisfied with $B$ and $E$ in $F$ and $A$, $C$ and $G$ in $Z_{2}$, $x$ and $%
m $ in $Z_{1}$, $l$ and $y$ in $H$ we have

$(1)$ $\eta (x,A)\le n^{-1}(p(m,G)+g(l,E)-g(y,B))$. For $C\in Z_2$ there
exist $P$ and $Q$ in $S_X$ with $\xi _0(0,Q)<0$ and $C\hat +P\le Q$. In view
of Condition $3.(iv)$ we have $B\hat + n_1\circ (a_1A)\hat +...\hat
+n_k\circ (a_kA) \hat +Q\le E\hat +G\hat +Q$. Then substituting $Q$ on $%
B\hat + n_1\circ (a_1A)\hat +...\hat +n_k\circ (a_kA) \hat +Q$ we get $B\hat
+n_1\circ (a_1A)\hat +...\hat +n_k\circ (a_kA) \hat +B\hat + n_1\circ
(a_1A)\hat +...\hat +n_k\circ (a_kA)\hat +Q=$ $2\circ B\hat +2\circ (
n_1\circ (a_1A)\hat +...\hat +n_k\circ (a_kA)) \hat +Q\le E\hat +G\hat
+Q\hat +B\hat + n_1\circ (a_1A)\hat +...\hat +n_k\circ (a_kA)$ $\le 2\circ
E\hat +2\circ G\hat +Q$ and by induction we get $k\circ B\hat +(k\circ (
n_1\circ (a_1A)\hat +...\hat +n_k\circ (a_kA)) \hat +Q) \le k\circ E \hat
+(k\circ G\hat +Q)$. Due to Lemma 2.5 for each $z\in Z_1$ there are the
following inequalities: $kn \eta (x,A)+\eta (z,Q)\le \eta (knx+z,k\circ
(n_1\circ (a_1A)\hat +...\hat +n_k\circ (a_kA)\hat +Q)$ $\le \xi
_0(knx+z,k\circ ( n_1\circ (a_1A)\hat +...\hat +n_k\circ (a_kA)\hat +Q)$ $%
\le \xi _0(km,k\circ G)+g(kl,k\circ E)-g(ky,k\circ B)$ $\le
k(p(m,G)+g(l,E)-g(y,B))+\xi _0(z,Q),$ so dividing by $kn$ and using that $%
k\in \mathbf{N}$ is arbitrary we get Inequality $(1)$.

$(iii)$. Necessity is evident, since $h\leq \xi _{0}$. We prove sufficiency.
In view of Theorem $2.8.(i)$ there exists a monotone Hahn-Banach extension $%
h $ of $g$ such that $\eta \leq h\leq \xi _{0}$, where $h$ is defined on $Z$
(see \S $2.8.(ii)$), since $\xi \leq \xi _{0}$ on $Z$. The rest of the proof
follows from Lemma 2.9 given below.

\textbf{2.9. Lemma.} \textit{If the preordering on} $X\times S$\textit{\ is
equality and there exists a superlinear map }$\eta :X\times S\rightarrow
(-\infty ,\infty ]$\textit{\ such that }$\eta \leq \xi _{0}$\textit{. A
finite semilinear and cancellative map }$g$\textit{\ defined on a }$\mathbf{P%
}$\textit{-quasisemilinear subspace }$E$\textit{\ of }$X\times S$\textit{\
has a finite semilinear extension }$h$\textit{\ to all of }$X\times S$%
\textit{.}

\textbf{Proof} is analogous to that of Lemma 2.2 with the cancellation
property instead of the monotone cancellation, that gives $h$ on $%
\bigtriangleup E$. Here is not demanded that $\eta \leq h\leq \xi _{0}$. If $%
X\times S\setminus \bigtriangleup E=:T\neq \emptyset ,$ then put $%
h(v,D)=\gamma $ for $(v,D)\in T$ and this gives $h$ on $sp_{\mathbf{K}%
}v\times Sp_{\mathbf{P}}D$, hence an extension $\bar{\xi}_{0}$ of $\xi _{0}$
from $T$ on $T\cup (sp_{\mathbf{K}}v\times Sp_{\mathbf{P}}D)$ is defined.
Choose $\gamma $ such that $\xi _{0}=\bar{\xi}_{0}$ on the extended in such
way a $\mathbf{P}$-quasisemilinear space $T_{v,D}$ generated by $T\cup
\{(v,D)\}$. Then, as above, we get $h$ on $\bigtriangleup T_{v,D}$.
Considering the family of all such extensions and applying to it the
Kuratowski-Zorn lemma as in \S 2.2 we get $h$ on $X\times S$.

\textbf{2.10. Theorem.} \textit{Let suppositions of Theorem 2.8 be satisfied
and }$X\times S$\textit{\ be monotonely cancellative and }$\eta \leq h\leq
\xi $\textit{. Suppose} $v\in X$\textit{\ and }$D\in S$\textit{\ and} $b\in
(-\infty ,\infty ]$\textit{. Then }$g$\textit{\ has a monotone Hahn-Banach
extension }$h$\textit{\ on }$X\times S$ \textit{with values in }$(-\infty
,\infty ]$\textit{\ such that }$h(v,D)=b$\textit{\ if and only if there
exists a superlinear map }$\eta :X\times S\rightarrow (-\infty ,\infty ]$ 
\textit{such that }

$(i)$ $\eta (x,A)\leq \xi (x+nv,A\hat{+}n_{1}\circ (a_{1}D)\hat{+}...\hat{+}%
n_{k}\circ (a_{k}D))-nb$\textit{\ for each }$x\in X$\textit{, }$A\in S$%
\textit{, }$k$\textit{\ and} $n_{i}$ \textit{in }$\mathbf{N}$\textit{, }$%
a_{i}\in \mathbf{K}_{+}$, $i=1,...,k$, $n=n_{1}a_{1}+...+n_{k}a_{k}$\textit{%
. If Condition }$(i)$\textit{\ is satisfied, then }$h$\textit{\ can be
chosen such that }$\eta \leq h\leq \xi $\textit{.}

\textbf{Proof.} Necessity is evident. We show sufficiency. Condition $(i)$
is the generalization of the condition $\xi \geq \eta $, since for $v=0$ and 
$D=\emptyset $ we get $\eta (x,A)\leq \xi (x,A)$. For the use of Lemma 2.6
in this situation we verify, that when $(v,D)\notin \bigtriangleup E$, then $%
U_{v,D}\geq b$. From $\eta \leq h\leq \xi $ on $E$ and Lemma $2.5$ we have $%
\xi (x_{2}+nv,A_{2}\hat{+}n_{1}\circ (a_{1}D)\hat{+}...\hat{+}n_{k}\circ
(a_{k}D))\leq $ $\xi (x_{2},A_{2})+n\xi (v,D)$ on $\bigtriangleup (E\cup sp_{%
\mathbf{K}}v\times Sp_{\mathbf{P}}D)$, consequently, $-\eta
(x_{2},A_{2})\geq -\xi (x_{2},A_{2})-n\xi (v,D)+nb,$ where $\xi (v,D)\geq b$%
. Then as in \S 2.6 we get $U_{v,D}\geq b$.

\textbf{2.11. Theorem.} \textit{Let }$X$\textit{, }$S$\textit{, }$p$\textit{,%
} $E$\textit{\ and }$g$\textit{\ be as in \S 2.8 and assume that }$X\times S$%
\textit{\ is monotonely cancellative and that }$(x,A)\geq (0,\emptyset )$%
\textit{\ for each }$x\in X$\textit{\ and }$A\in S$\textit{. Then }$g$%
\textit{\ has a monotone Hahn-Banach extension if and only if }

$(i)$ $g(x,A)\leq g(y,B)+p(z,C)$\textit{\ whenever }$(x,A)\leq (y+z,B\hat{+}%
C)$\textit{\ with }$(x,A)$\textit{\ and }$(y,B)$\textit{\ in }$E$\textit{$.$
If }$p$\textit{\ is monotone, then for each }$v\in X$\textit{\ and }$D\in S$ 
\textit{there exists a monotone and additive map }$h$\textit{\ on }$X\times
S $\textit{\ such that} $h\leq p$\textit{\ and }$h(v,D)=p(v,D)$\textit{.}

\textbf{Proof.} Condition $2.10.(i)$ for $b=\xi (v,D)<\infty $ can be
written in the form $\eta (x,A)\le \lim_{n\to \infty } (\xi (x+nv,A\hat
+n_1\circ (a_1D)\hat +...\hat +n_k\circ (a_kD))- n\xi (v,D))$, where $x\in X$
and $A\in S$.

\textbf{2.12. Corollary.} \textit{Let }$X\times S$\textit{\ be a monotonely
cancellative }$\mathbf{P}$\textit{-quasisemilinear preordered space and }$p$%
\textit{, }$\xi $\textit{, }$g$\textit{\ be as usual. Suppose that for a
given }$(v,D)\in X\times S$\textit{\ }

$(i)$ $\inf \{\xi (x+y+v,A\hat{+}B\hat{+}D)-\xi (y+v,B\hat{+}D):$ $(y,B)\in
X\times S\}>-\infty $\textit{\ for each }$(x,A)\in X\times S$\textit{\ and }

$(ii)$ $-\infty <\xi (y+v,B\hat{+}D)<\infty $ \textit{for some }$(y,B)\in
X\times S$\textit{. Then there exists a monotone Hahn-Banach extension} $h$%
\textit{\ of }$g$\textit{\ such that} $h(v,D)=\xi (v,D)$\textit{.}

\textbf{Proof.} For $(x,A)\in X\times S$ put $\eta (x,A)=\inf \{ \xi
(x+y+v,A\hat +B\hat +D)$ $-\xi (y+v,B\hat +D):$ $(y,B)\in X\times S \} $. By
Condition $(i)$ $\eta >-\infty $. Then as in \cite{topsoe} it can be shown
that $\eta \le \xi $ and $\eta $ is superlinear, since $\xi $ is sublinear
and due to Theorem 2.10 we get the statement of this Corollary.

\textbf{2.13. Definitions and Notes.} Let $X\times S$ be a $\mathbf{P}$%
-quasisemilinear preordered space, $p: X\times S\to (-\infty ,\infty ]$ be
sublinear, $E$ be a $\mathbf{P}$-quasisemilinear subspace and $g: E\to 
\mathbf{R}$ be semilinear on $E$. Consider a subset $W$ of $X\times S$. Then 
$h: X\times S\to (-\infty , \infty ]$ is called a $W$-maximal monotone
Hahn-Banach extension of $g$ if for each monotone Hahn-Banach extension $%
h^{\prime}$ of $g$ from the inequality $h^{\prime}\ge h$ on $W$ it follows $%
h^{\prime}=h$ on $W$.

A sublinear map $\xi :X\times S\rightarrow (-\infty ,\infty ]$ is called of
moderate variation if%
$\inf \{\xi (x+y,A\hat{+}B)-\xi (y,B):(y,B)\in X\times S\}>-\infty $
for all $(x,A)\in X\times S$. For $\xi $ of moderate variation define a
superlinear map $\eta $ by%
$\eta (x,A)=\inf \{\xi (x+y,A\hat{+}B)-\xi (y,B):(y,B)\in X\times S\},$
where $(x,A)\in X\times S$.

\textbf{2.14. Theorem.} \textit{Let }$X$\textit{, }$S$\textit{,} $p$\textit{%
, }$E$\textit{\ be as in \S 2.13 and assume that }$X\times S$\textit{\ is
monotonely cancellative. If }$\xi =\xi _{p,g}$\textit{\ is a sublinear map
of moderate variation, then for each subset }$W$\textit{\ of }$X\times S$%
\textit{\ there exists a }$W$\textit{-maximal monotone Hahn-Banach extension
of }$g$\textit{.}

\textbf{Proof.} Consider $\eta (x+y,A\hat +B)= \inf \{ \xi (x+y+v,A\hat
+B\hat +D)-$ $\xi (v,D):$ $(v,D) \in X\times S \} $, where $(x,A)$ and $%
(y,B) $ are in $X\times S$. We have $\eta (x,A)+\eta (y,B)=\inf \{ \xi
(x+v_1,A\hat +D_1)+$ $\xi (y+v_2,B\hat +D_2)-$ $\xi (v_1,D_1)-$ $\xi
(v_2,D_2):$ $(v_1,D_1)$ $\mbox{and}$ $(v_2,D_2)\in X\times S \} \le $ $\inf
\{ \xi (x+y+(v_1+v_2),A\hat +B\hat +(D_1\hat +D_2))-$ $\xi (v_1,D_1)-$ $\xi
(v_2,D_2):$ $(v_1,D_1)$ $\mbox{and}$ $(v_2,D_2)\in X\times S \} \le $ $\inf
\{ \xi (x+y+v_1+0,A\hat +B\hat +D_1\hat +\emptyset )-$ $\xi (v_1,D_1)-$ $\xi
(0,\emptyset )$ $(v_1,D_1)\in X\times S \} $ $=\eta (x+y,A\hat +B)$, hence $%
\eta $ satisfies $2.1.2.(S3)^{\prime},(S4)^{\prime}$. Since $\xi $ satisfies 
$2.1.2.(S6)$, then $\eta $ also satisfies $(S6)$. Since $\xi $ satisfies $%
(S2)$, then $\eta $ also satisfies $(S2)$. Evidently, $\eta (0,\emptyset )=0$%
, since $\xi (0,\emptyset )=0$. From $\xi (x+y,A\hat +B)\le \xi (x,A)+\xi
(y,B)$ it follows $\eta (x,A)\le \xi (x+y,A\hat +B)-\xi (y,B)\le \xi (x,A)$.
Taking $p=\xi $, $E=\{ 0,\emptyset \} $ and $b=\xi (v,D)$ we infer from
Theorem 2.10 that for each $v\in X$ and each $D\in S$ there exists a
monotone semilinear map $h: X\times S\to (-\infty ,\infty ]$ such that $\eta
\le h\le \xi $ and $h(v,D)=\xi (v,D)$. The map $\xi $ is of moderate
variation, hence $V:=\{ (x,A):$ $\xi (x,A)<\infty ,$ $(x,A)\in X\times S \}
= $ $\{ (x,A):$ $\eta (x,A)<\infty ,$ $(x,A)\in X\times S \} $ and $V$ is a
hereditary face, since $\xi $ is monotonely cancellative in accordance with
Lemma $2.5.(i)$.

Suppose $\xi <\infty $ and $(v,D)\in X\times S$, consider the set $F:= \{
(x+y,A\hat +B)):$ $(x,A)\in E,$ $x\in sp_{\mathbf{K}}v,$ $B\in Sp_{\mathbf{P}%
}D \} $. Define $h$ on $F$ by $h(x+nv,A\hat +n_1\circ (a_1D)\hat +...\hat
+n_k\circ (a_kD))=$ $g(x,A)+n\xi (v,D)$ for each $k$ and $n_i\in \mathbf{N_0}
$ and $a_i\in \mathbf{K_+}$ and each $(x,A)\in E$, where $%
n=n_1a_1+...+n_ka_k $. Put $\xi ^{\prime}:=\xi _{p,h}$, then $\xi ^{\prime}$
is of moderate variation and the associated superlinear map dominates $\eta $%
. Moreover, if $r: X\times S\to (-\infty ,\infty ]$ is a monotone
Hahn-Banach extension of $g$ and if $r\ge h$ on $F\cap W$, then $r=h$ on $%
F\cap W$. From the definition of $\xi ^{\prime}$ it follows that

$\xi ^{\prime}(x,A)=\min \{ \xi (x,A), C_1, C_2 \} $, where

$C_1:=\inf \{ m^{-1}(p(z,Q)+g(l,J)-g(y,B) +n\xi (v,D)):$ $(y+mx,B\hat
+m_1\circ (a_1A)\hat +...\hat +m_j\circ (a_jA))\le $ $(l+nv+z,J\hat
+n_1\circ (d_1D)\hat +...\hat +n_k\circ (d_kD)\hat +Q) \} ,$

$C_2:=\inf \{ m^{-1}(p(z,Q)+g(l,J)-g(y,B) -n\xi (v,D)):$ $(y+nv+mx,B\hat
+m_1\circ (a_1A)\hat +...\hat +m_j\circ (a_jA)$ $\hat +n_1\circ (d_1D)\hat
+...\hat +n_k\circ (d_kD))\le $ $(l+z,J\hat +Q) \} ,$ where $%
m=m_1a_1+...+m_ja_j$, $n=n_1d_1+...n_kd_k$; $j$, $k$, $m_i$ and $n_i\in 
\mathbf{N}$; $a_i$ and $d_i \in \mathbf{K_+}$; $(y,B)$ and $(l,J)\in E$, $%
m\in \mathbf{N}$. In view of Lemma $2.5.(v)$ $C_1\ge \xi (v,D)$, since $%
p(z,Q)+n\xi (v,D)\ge \xi (z+nv,Q\hat +n_1\circ (d_1D) \hat +...\hat
+n_k\circ (d_kD)).$ Analogously to Lemma 10 \cite{topsoe} we get:

$(i)$ $\xi ^{\prime}(x,A)=\inf \{ m^{-1}(\xi (mx+nv,$ $m_1\circ (a_1A)\hat
+...\hat +m_j\circ (a_jA))\hat +$ $n_1\circ (d_1D)\hat +...\hat +n_k\circ
(d_kD))-$ $\xi (nv,n_1\circ (d_1D)\hat +...\hat +n_k\circ (d_kD)):$ $j, k,
m_i, n_i, a_i, d_i \} $. \newline
Thus Equation $(i)$ can be rewritten as:

$(ii)$ $\xi ^{\prime }(x,A)=\lim_{n/m\rightarrow \infty }\{m^{-1}(\xi
(mx+nv, $ $m_{1}\circ (a_{1}A)\hat{+}...\hat{+}m_{j}\circ (a_{j}A))\hat{+}$ $%
n_{1}\circ (d_{1}D)\hat{+}...$ $\hat{+}n_{k}\circ (d_{k}D))-$ $\xi
(nv,n_{1}\circ (d_{1}D)\hat{+}...\hat{+}n_{k}\circ (d_{k}D)):$ $%
j,k,m_{i},n_{i},a_{i},d_{i}\}$ with natural numbers $m$ and $n$, since $\xi $
satisfies $2.1.2.(S6)$, hence

$\xi ^{\prime}(x+z,A\hat +Q)-\xi ^{\prime}(z,Q)=$ $\lim_{n/m\to \infty } \{
m^{-1}(\xi (mx+mz+nv,$ $m_1\circ (a_1A)\hat +...\hat +m_j\circ (a_jA))\hat
+m\circ Q\hat +$ $n_1\circ (d_1D)\hat +...\hat +n_k\circ (d_kD))-$ $\xi
(mz+nv,m\circ Q\hat + n_1\circ (d_1D)\hat +...\hat +n_k\circ (d_kD)):$ $j,
k, m_i, n_i, a_i, d_i \} $, and inevitably $\xi ^{\prime}$ is of moderate
variation and the associated superlinear map dominates $\eta $. The rest of
the proof of Theorem 2.14 is analogous to that of Theorem 5 \cite{topsoe}
with the help of lemmas given above.

Addresses:

S.V. L\"{u}dkovsky

Theoretical Department, Institute of General Physics, Russian Academy of
Sciences,

Str. Vavilov 38, Moscow, 119991 GSP-1, Russia.

E-mail address: ludkovsk@fpl.gpi.ru \newline

J.C. Ferrando

Centro de Investigaci\'{o}n Operativa, Universidad Miguel Hern\'{a}ndez,

E-03202 Elche (Alicante), Comunidad Valenciana, Spain.

E-mail address: jc.ferrando@umh.es


\begin{thebibliography}{99}
\bibitem{att} Attouch, H., Lucchetti, R., Wets, R.J.-B. ''The topology of
the $\rho $-Hausdorff distance''//Annali di Mat. Pura ed Appl. Ser IV. 1991.
V.\textbf{160}, P. 303-320.

\bibitem{ferr} Ferrando, J.C., L\'{o}pez Pellicer, M., S\'{a}nchez Ruiz,
L.M. ''Metrizable barrelled spaces''. Pitman Research Notes in Math. V. 
\textbf{332} (Longman: Harlow, Essex, UK, 1995).

\bibitem{luumn48} L\"{u}dkovsky, S.V. ''Topological groups and their $\kappa 
$-metrics''// Usp. Mat. Nauk. 1993. V. \textbf{48}, N 1, 173-174.

\bibitem{lusmz} L\"{u}dkovsky, S.V. ''$\kappa $-normed topological vector
spaces''// Sibirsk. Mat. J. 2000. V. \textbf{41}, N 1, 167-184.

\bibitem{ludkns} L\"{u}dkovsky, S.V. ''Duality of $\kappa $-normed
topological vector spaces and their applications'', Los Alamos National
Laboratory, USA. Preprint \textbf{math.GN/0101088}, 26 pages, 10 January
2001 (shortly: November 2000, abstract of P.S. Alexandroff's seminar, Vestn.
Mosk. State Univ. Ser. Mat.).

\bibitem{nari} Narici, L., Beckenstein, E. ''Topological vector spaces''.
(Marcel-Dekker Inc.: New York, 1985).

\bibitem{schaef} Schaefer, H.H., Wolff, M.P. ''Topological vector spaces''
(Springer: New York, Sec. Edit., 1999).

\bibitem{shep1} Shepin, E.V. ''Topology of the limit spaces of uncountable
inverse spectra''// Usp. Math. Nauk. 1976. V. \textbf{31}, N 5, P. 191-226.

\bibitem{shep2} Shepin, E.V. ''About $\kappa $-metrizable spaces''// Izv.
Akad. Nauk SSSR. Ser. Math. 1979. V. \textbf{43}, N 2, P. 442-477.

\bibitem{topsoe} Topsoe, F. ''The naive approach to the Hahn-Banach
theorem''// Commentationes Math. 1978. V. \textbf{21}, P. 315-330.
\end{thebibliography}
\end{document}